\documentclass{fpsac}

\newtheorem{theorem}{Theorem}
\newtheorem{lemma}[theorem]{Lemma}

\newtheorem{conj}[theorem]{Conjecture}

\def\mT{\mathcal{T}}

\title{Proof of the Somos-4 Hankel Determinants Conjecture}

\author[G. Xin]{Guoce Xin}
\address{Center for
Combinatorics, LPMC-TJKLC, Nankai University, Tianjin 300071, P. R.
China} \email{gxin@nankai.edu.cn}



\begin{document}

\begin{abstract}
By considering the fundamental equation $x=y-y^2=z-z^3$, Somos
conjectured that the Hankel determinants for the generating series
$y(z)$ are the Somos-4 numbers. We prove this conjecture by using
the quadratic transformation for Hankel determinants of Sulanke and
Xin.
\end{abstract}

\maketitle  

\section{Introduction}

A generating function $Q(x)=\sum_{n\ge 0} q_n x^n$ defines a
sequence of Hankel matrices $H_1,H_2,H_3,\dots, $ where $H_n$ is an
$n$ by $n$ matrix with entries $(H_n)_{i,j}=q_{i+j-2}$. Hankel
determinants are determinants of these matrices. Traditionally,
$H_0$ is defined to be the empty matrix with determinant $1$.

In the year of 2000, Somos \cite{Somos} considered the fundamental
equation $x=y-y^2=z-z^3$. He observed that three types of expansions
give nice Hankel determinants. The first one is by expanding $y$ as
a series in $x$, which gives the generating function for Catalan
numbers; the second one is by expanding $y$ as a series in $z$,
which gives a generating function related to Catalan and Motzkin
numbers; the third one is by expanding $z$ as a series in $x$, which
gives the generating function for ternary trees. The first case was
known by Shapiro \cite{Shapiro}, the third case was proved
independently in \cite{err,xingessel2,tamm}, and the second case,
known as the \emph{Somos-4 conjecture}, is still open.

The Somos-4 conjecture can be restated as follows. Expanding $y$ as
a series in $z$ gives
$$y=z+z^2+z^3+3z^4+8z^5+23z^6 +\cdots.$$
Let $Q(z)=(y-z)/z^2$ and let $s_n=\det H_n(Q)$.

\begin{conj}[Somos-4]
The Hankel determinants $s_n$ defined above satisfy the recursion
\begin{align} s_ns_{n-4}=s_{n-1}s_{n-3}+s_{n-2}^2, \label{e-rec-s}
\end{align} with
initial conditions $s_0=1,s_1=1,s_2=2,s_3=3$.
\end{conj}

For instance,
$$H_3(Q)= \left(
                 \begin{array}{ccc}
                   1 & 1 & 3 \\
                   1 & 3 & 8 \\
                   3 & 8 & 23 \\
                 \end{array}
               \right), \qquad s_3=\det H_3(Q)=3.
$$
Our main objective in this paper is to prove the above conjecture.

There are many classical tools of continued fractions for evaluating
Hankel determinants, such as the $J$-fractions in Krattenthaler
\cite{kratt} or Wall \cite{wall} and the $S$-fractions in Jones and
Thron \cite[Theorem~7.2]{cfraction}. Our tool is by  Sulanke and
Xin's quadratic transformation for Hankel determinants
\cite{Sulanke-xin} developed from the continued fraction method of
Gessel and Xin \cite{xingessel2}.

\section{Solving a system of recurrences}
\def\mT{\mathcal{T}}

Proposition 4.1 of \cite{Sulanke-xin} defines a quadratic
transformation $\mT$, and asserts that for certain generating
function $F$, we can find $\mT(F)$ such that $\det(H_n(F))=a \det
(H_{n-d-1}(\mT(F)))$, where $a$ is a constant and $d$ is a
nonnegative integer. See \cite{Sulanke-xin} for detailed
information. Here we only need the following special case.

\begin{lemma}\label{l-transform}
Suppose $a\ne 0$. If the generating functions $F(x)$ and $G(x)$ are
uniquely defined by \begin{align*}
 F(x)&=
\frac{a+bx}{1+cx+dx^2+x^2(e+fx)F(x)},\\
G(x)&=\frac{-{\frac {{a}^{3}e+{a}^{2}d-acb+{b}^{2}}{{a}^{2}}}-{\frac
{{a}^{4}f+c
{a}^{3}d-{c}^{2}{a}^{2}b+2\,ca{b}^{2}-b{a}^{2}d-{b}^{3}}{{a}^{3}}}x}{1+cx-
{\frac {-2\,acb+2\,{b}^{2}+{a}^{2}d}{{a}^{2}}}x^2+x^2(-1-{\frac
{b}{a}}x)G(x)} ,
\end{align*}
then $\det H_n(F) =a^n \det H_{n-1}(G)$.
\end{lemma}

Our proof is by iterative application the above lemma. To be
precise, define $Q_0(x)=Q(x)$, and recursively define $Q_{n+1}(x)$
to be the unique power series solution of
\begin{align}
Q_{n+1}(x)
=\frac{a_{n+1}+b_{n+1}x}{1+c_{n+1}x+d_{n+1}x^2+x^2(e_{n+1}+f_{n+1}x)Q_{n+1}(x)},
\end{align}
where
\begin{align}
a_{n+1}&=-{\frac {{a_{n}}^{3}e_{n}+{a_{n}}^{2}d_{n}-a_{n}c_{n}b_{n}+{b_{n}}^{2}}{{a_{n}}^{2}}} \label{e-a}\\
b_{n+1}&=-{\frac {{a_{n}}^{4}f_{n}+c_{n}
{a_{n}}^{3}d_{n}-{c_{n}}^{2}{a_{n}}^{2}b_{n}+2\,c_{n}a_{n}{b_{n}}^{2}-b_{n}{a_{n}}^{2}d_{n}-
{b_{n}}^{3}}{{a_{n}}^{3}}}\label{e-b}
\\
c_{n+1}&=c_{n}\nonumber \\
d_{n+1}&=- {\frac {-2\,a_{n}c_{n}b_{n}+2\,{b_{n}}^{2}+{a_{n}}^{2}d_{n}}{{a_{n}}^{2}}}\label{e-d}\\
e_{n+1}&=-1\nonumber \\
f_{n+1}&=-{\frac {b_{n}}{a_{n}}}\label{e-f}
\end{align}

It is straightforward to represent $Q(x)$ as the unique power series
solution of
$$ Q(x)= \frac{1-x}{1-2x-x^2Q(x)}. $$
Therefore we shall set $a_0=1,b_0=-1,c_0=-2,d_0=0,e_0=-1,f_0=0.$ By
Lemma \ref{l-transform}, one can deduce that
$\det(H_n(Q))=a_0^na_1^{n-1}\cdots a_{n-1}$. This transforms the
recursion for $s_n$ to that for $a_n$ as follows:
\begin{align}
\label{e-rec-a} a_na_{n-1}a_{n-2}=1+1/a_{n-1}.
\end{align}
We remark that the above recursion implies that $s_3=s_2+s_1^2$,
which holds for the Somos-4 sequence.

It is a surprise that the recursion system can be solved for
arbitrary initial condition. For simplicity, we write $c_n=c$ and
assume $e_0=-1$ (otherwise start with $Q_1$). Our solution can be
stated as follows.
\begin{theorem}
\label{t-rec-a} Suppose $c_n=c$, $e_n=-1$, and $a_n,b_n,d_n,f_n$
satisfy the recursion (\ref{e-a},\ref{e-b},\ref{e-d},\ref{e-f}).
Then 
\begin{align}
\label{e-t-rec-a} a_{n+2}a_{n+1}+a_{n+1}a_n =2a_0a_1
+a_0(f_0+f_1+c)(2f_1+c)-(a_0(f_0+f_1+c))^2/a_{n+1}.
\end{align}
\end{theorem}

\begin{proof}
We shall try to write everything in terms of the $a$'s. Using
\eqref{e-f}, we can replace $b_n$ with $-a_n f_{n+1}$ everywhere.
Therefore \eqref{e-a} becomes
\begin{align}
\label{e-d-new} d_n=a_n-a_{n+1}-cf_{n+1}-f_{n+1}^2.
\end{align}
Substituting \eqref{e-d-new} into \eqref{e-b} and simplifying gives
$$f_{n+2}a_{n+1}=a_nf_n+ca_n-ca_{n+1}+f_{n+1}a_n-f_{n+1}a_{n+1}, $$
which can be written as
$$a_{n+1}(f_{n+2}+f_{n+1}+c)=a_n(f_{n+1}+f_n+c). $$
That is to say
\begin{align}
\label{e-ffc} a_{n+1}(f_{n+2}+f_{n+1}+c)=a_0(f_1+f_0+c).
\end{align}

Substituting \eqref{e-d-new} into \eqref{e-d} and simplifying gives
$$a_n-a_{n+2}=cf_{n+2}+f_{n+2}^2 -(cf_{n+1}+f_{n+1}^2)=(f_{n+2}-f_{n+1})(f_{n+2}+f_{n+1}+c).$$
Applying \eqref{e-ffc}, we obtain
$$
a_na_{n+1}-a_{n+1}a_{n+2} =a_0(f_1+f_0+c) (f_{n+2}-f_{n+1}),
$$
which leads to
\begin{align}\label{e-fff}
a_0a_{1}-a_{n+1}a_{n+2} =a_0(f_1+f_0+c) (f_{n+2}-f_{1}).
\end{align}
Combining \eqref{e-ffc} and \eqref{e-fff}, we obtain
\eqref{e-t-rec-a}.
\end{proof}

Now we are ready to prove the Somos-4 Conjecture.
\begin{proof}[Proof of the Somos-4 Conjecture]
Applying Theorem \ref{t-rec-a} for the case
$a_0=1,b_0=-1,c=-2,d_0=0,e_0=-1,f_0=0$, we obtain $a_1=2,$ $ f_1=1$,
and
\begin{align}\label{e-an2}
a_{n+2}=4/a_{n+1}-a_n-1/a_{n+1}^2.
\end{align}

Recall that we have transformed the recursion \eqref{e-rec-s} to
\eqref{e-rec-a}, which can be written as
$$a_na_{n-1}^2a_{n-2}-1-a_{n-1}=0. $$
By applying \eqref{e-an2} (with $n$ replaced by $n-2$) and
simplifying, the above equation becomes
$$4a_{n-2}a_{n-1}-a_{n-2}-a_{n-2}^2a_{n-1}^2-1-a_{n-1}=0. $$
Denote by $T(n)$ the left-hand side of the above equation. We claim
that $T(n)=0$ for all $n$, so that \eqref{e-rec-a} holds and the
conjecture follows.

We prove the claim by induction on $n$. The claim is easily checked
to be true for $n=2$. Assume the claim hold for $n-1$. By applying
\eqref{e-an2} (with $n$ replaced by $n-3$) and simplifying, we
obtain
$$T(n)=4a_{n-3}a_{n-2}-a_{n-2}-a_{n-3}-a_{n-3}^2a_{n-2}^2-1=T(n-1)=0. $$
Thus the claim follows.
\end{proof}


\noindent {\bf Acknowledgments.} The author was grateful to Doron
Zeilberger for calling his attention to the Somos-4 conjecture. This
work was supported by the 973 Project, the PCSIRT project of the
Ministry of Education, the Ministry of Science and Technology and
the National Science Foundation of China.

\bibliographystyle{amsplain}

\begin{thebibliography}{10}
\renewcommand\emph[1]{{#1}}

\bibitem{err}
{\"O}. E{\u{g}}ecio{\u{g}}lu, T. Redmond, and C. Ryavec, \emph{From
a polynomial {R}iemann hypothesis to alternating sign matrices},
Electron. J. Combin. \textbf{8} (2001), no.~1,  R36, 51 pp.


\bibitem{xingessel2}
I.~M.~Gessel and G.~Xin, \emph{The generating function of ternary
trees and
  continued fractions}, Electron. J. Combin., \textbf{13} (2006), R53. (electronic).

\bibitem{cfraction}
W.~B.~Jones and W.~J.~Thron, \emph{Continued Fractions: Analytic
Theory and Applications},
  Encyclopedia of Mathematics and its Applications, vol.~11, Addison-Wesley, Reading, Mass., 1980.

\bibitem{kratt}
 C.~Krattenthaler,
 Advanced determinant calculus: a complement, \emph{Linear Algebra
     Appl.} \textbf{411} (2005), 68--166

The Andrews Festschrift (Maratea, 1998). \emph{Sem. Lothar. Combin.}
\textbf{42} (1999), Art. B42q, 67 pp. (electronic).


\bibitem{Shapiro}  L.~W.~Shapiro,  A Catalan triangle. \emph{Discrete Math.} \textbf{14} no. 1  (1976),  83--90.

\bibitem{Somos}
M.~Somos, {\tt
http://grail.cba.csuohio.edu/\raisebox{-3pt}{\~{}}somos/nwic.html}.

\bibitem{Sulanke-xin}  R.~A.~Sulanke and G.~Xin,
Hankel Determinants for Some Common Lattice Paths, \emph{Adv. in
Appl. Math.}, to appear, appeared at Formal Power Series and
Algebraic Combinatorics (FPSAC06).

\bibitem{tamm}
U. Tamm, \emph{Some aspects of {H}ankel matrices in coding theory
and combinatorics}, Electron. J. Combin. \textbf{8} (2001), no.~1,
A1, 31 pp.


\bibitem{wall}
H.~S.~Wall, \emph{Analytic Theory of Continued Fractions}, Van
Nostrand, New York, 1948.


\end{thebibliography}

\providecommand{\bysame}{\leavevmode\hbox
to3em{\hrulefill}\thinspace}
\providecommand{\MR}{\relax\ifhmode\unskip\space\fi MR }
\providecommand{\MRhref}[2]{%
  \href{http://www.ams.org/mathscinet-getitem?mr=#1}{#2}
} \providecommand{\href}[2]{#2}

\end{document}